\documentclass[12pt]{amsart}
\usepackage{amsfonts,latexsym}

\newtheorem{theorem}{Theorem}

\newtheorem{proposition}[theorem]{Proposition} 
 
\newtheorem{corollary}[theorem]{Corollary}

\theoremstyle{definition}
\newtheorem{definition}[theorem]{Definition}

\theoremstyle{remark}

\DeclareMathOperator{\cf}{cf}
\newcommand{\cupdoty}{\dot{\cup}} 
\newcommand{\bigcupdoty}{\dot{\bigcup}} 

\title[Combinatorial principles, compactness of spaces III]
{Combinatorial and model-theoretical principles related to 
regularity of ultrafilters and compactness of topological spaces. III.}

\author[]{Paolo Lipparini} 
\address{Dipartimento di Matematica\\
Viale della Ricerca Scientifica\\
II Universitaccia Romanaccia (Tor Vergata)\\
I-00133 ROME 
ITALY
}
\urladdr{http://www.mat.uniroma2.it/\textasciitilde lipparin}
\thanks{The author has received support from MPI and GNSAGA.
We wish to express our gratitude to X. Caicedo for stimulating discussions and correspondence} 


\begin{document} 

\begin{abstract} 
We generalize the results from \cite{topproc}; in particular
the present results apply to singular cardinals, too.
\end{abstract}

\maketitle




See \cite{parti,KV,EGT} for definitions and notation. 

We shall need the following theorem from \cite{nd}. 

\begin{theorem} \label{limit} 
If $ \lambda $ is a singular cardinal, then 
an ultrafilter is $(\lambda ,\lambda )$-regular
if and only if it is either 
$(\cf\lambda ,\cf\lambda )$-regular
or
$(\lambda^+ ,\lambda^+ )$-regular. 
\end{theorem}

\begin{corollary} \label{cor}
Suppose that $ \lambda $ is a singular cardinal,
and consider the topological space $X$ 
obtained by forming the disjoint union of the topological 
spaces $ \lambda ^+$ and $\cf \lambda $, both
endowed with the order topology. 

Then, for every ultrafilter $D$, the space  $X$
is $D$-compact if and only if $D$ is not 
$( \lambda, \lambda )$-regular. 
 
Thus, $X$ is productively $[ \lambda', \mu']$-compact if and only if
there exists a 
$( \lambda', \mu' )$-regular
 not $( \lambda, \lambda )$-regular
ultrafilter. In particular, $X$ is not productively $[ \lambda, \lambda ]$-compact. 
 \end{corollary} 

\begin{proof}
By Theorem \ref{limit}, $D$ is not 
$( \lambda, \lambda )$-regular
if and only if it is neither 
$(\cf\lambda ,\cf\lambda )$-regular
nor $(\lambda^+ ,\lambda^+ )$-regular. 

Hence, by \cite[Proposition 1]{topproc}, and since both
$ \lambda ^+$ and $\cf \lambda $ are regular cardinals,
$D$ is not 
$( \lambda, \lambda )$-regular
if and only if both $ \lambda ^+$
and $\cf \lambda $ are $D$-compact.
 This is clearly equivalent to $X$
being $D$-compact.

The last statement is immediate from \cite[Theorem 1.7]{C}, also stated
in \cite[Theorem 2]{topproc}.  
 \end{proof} 

Let $\mathbf{2}= \{0,1 \} $ denote the two-elements topological space
with the discrete topology. 
If $ \lambda \leq \mu$ are cardinals, let 
$\mathbf{2}^ \mu $
be the Tychonoff product of $ \mu $-many copies
of $\mathbf{2}$, and let
$\mathbf{2}^ \mu_ \lambda  $
denote the subset of
$\mathbf{2}^ \mu $
consisting of all those functions
$h: \mu \to  \mathbf{2}$ such that
$ \left| \{ \alpha \in \mu| h( \alpha )=1 \} \right| < \lambda $.

In passing, let us mention that, when $ \mu= \aleph_\omega$,
the space $\mathbf{2}^ \mu_ \mu  $ provides
an example of a linearly Lindel\"of not Lindel\"of space.
See \cite[Example 4.1]{AB}. Compare also \cite[Example 4.2]{St}.

Notice that 
$\mathbf{2}^ \mu_ \lambda  $ is a Tychonoff topological group
with a base of clopen sets.

Set theoretically, $\mathbf{2}^ \mu_ \lambda  $
is in a one to one correspondence (via characteristic functions) with
$S_\lambda(\mu)$, the set of all subsets of $ \mu$
of cardinality $ < \lambda $. 
Since many properties of ultrafilters  are defined in terms of $S_\lambda(\mu)$,
for sake of convenience, in what follows we shall deal with
$S_\lambda(\mu)$, rather than $\mathbf{2}^ \mu_ \lambda  $.
Henceforth, we shall deal with the topology induced on $S_\lambda(\mu)$
by the above correspondence.

In detail, $S_\lambda(\mu)$ is endowed with the smallest topology containing, as open sets, all sets of the form $X_\alpha = \{x\in S_\lambda(\mu) | \alpha\in x\}$ ($\alpha$ varying in $\mu$), as well as their complements. Thus, a base for the topology consists of all finite intersections of the above sets; that is, the elements of the base are the sets $\{x\in S_\lambda(\mu) | \alpha_1\in x,  \alpha_2\in x , \dots, \alpha_n\in x , \beta_1\not\in x,  \beta_2\not\in x , \dots, \beta_m\not\in x \}$,  with  $n,m$ varying in $\omega$ and $\alpha_1, \dots, \alpha_n, \beta_1 , \dots, \beta_m
$ varying $\mu$. 

Notice that this topology is finer than the topology on $S_\lambda(\mu)$ used in
\cite{topproc}.

With the above topology, $S_\lambda(\mu)$ and
$\mathbf{2}^ \mu_ \lambda  $ are homeomorphic,
thus $S_\lambda(\mu)$ can be given the structure 
of a Tychonoff topological group.

Notice that if $ \lambda \leq \mu$ then $S_\lambda(\mu)$
is not $[\lambda ,\lambda ]$-compact. Indeed, for $ \alpha \in \mu$,
let $Y_\alpha = \{ x \in S_\lambda(\mu)| \alpha \not\in x \} $.
If $Z \subseteq \mu$ and $|Z|= \lambda $ then
$(Y_\alpha )_{\alpha\in Z}$ is an open cover of 
$S_\lambda(\mu)$ by $\lambda $-many sets,
$<\lambda $ of which never cover $S_\lambda(\mu)$. 

\begin{proposition} \label{dsll}  
For every  ultrafilter $D$ and every cardinal $\lambda$, 
the topological space $S_\lambda(\lambda)$ is $D$-compact if and only if 
 $D$ is not ($\lambda,\lambda$)-regular.
\end{proposition}

\begin{proof}  
Suppose that  
$D$ is an ultrafilter over $I$ and that
$S_\lambda(\lambda)$ is $D$-compact. For every $f:I\to S_\lambda(\lambda)$ there exists $x\in S_\lambda(\lambda)$ such that $f(i)_{i\in I}$ $D$-converges to $x$. If $\alpha\in\lambda  $ and $\{i\in I|\alpha\in f(i)\}\in D$ then $\alpha\in x $, since 
 otherwise $Y=\{z\in S_\lambda(\lambda) | \alpha\not\in z\}$ is an open set containing $x$, and $\{i\in I|f(i)\in Y\}= \{i\in I|\alpha\not\in f(i)\}\not\in D$, contradicting $D$-convergence.

Whence, $\{\alpha\in\lambda  |\{i\in I|\alpha\in f(i)\}\in D\}\subseteq x\in S_\lambda(\lambda)$, and thus $x$ has cardinality $<\lambda$; that is, $f$
does not witness ($\lambda,\lambda$)-regularity of $D$.
Since $f$ has been chosen arbitrarily,
$D$ is not ($\lambda,\lambda$)-regular.

Conversely, suppose that $D$ over $I$ is not ($\lambda,\lambda$)-regular, and let 
$f:I\to S_\lambda(\lambda)$. Then $x=\{\alpha\in\lambda  |\{i\in I|\alpha\in f(i)\}\in D\}$ has cardinality $<\lambda$ and hence is in $S_\lambda(\lambda)$. We show that $f$ $D$-converges to $x$. Indeed, let $Y$ be a neighborhood of $x$: we have to show that 
 $\{i\in I|f(i)\in Y\}\in D$. 
Without loss of generality, we can suppose that $Y$ is an element of the base 
of $S_\lambda(\lambda)$, that is, $Y$ has the
 form $\{z\in S_\lambda(\lambda) | \alpha_1\in z,  \alpha_2\in z , \dots, \alpha_n\in z , \beta_1\not\in z,  \beta_2\not\in z , \dots, \beta_m\not\in z \}$. Since $D$ is closed under finite intersections,
 then $\{i\in I|f(i)\in Y\}\in D$ if and only if    $\{i\in I|\alpha_1\in f(i)\}\in D$ and
 $\{i\in I|\alpha_2\in f(i)\}\in D$ and\dots\ and  $\{i\in I|\alpha_n\in f(i)\}\in D$ and 
 $\{i\in I|\beta_1\not\in f(i)\}\in D$ and\dots\ and $\{i\in I|\beta_m\not\in f(i)\}\in D$. But all 
the above sets are actualy in $D$, by the definition of $x$ and since $x\in Y$ and 
$D$ is an ultrafilter; thus $f$ $D$-converges to $x$.

Since $f$ was arbitrary, every $f:I\to S_\lambda(\lambda)$ $D$-converges, and thus 
$S_\lambda(\lambda)$ is $D$-compact. 
\end{proof}

\begin{corollary} \label{psll}  
 The space $S_\lambda(\lambda)$ is productively $[\lambda',\mu']$-compact if and only if  there exists a $(\lambda',\mu')$-regular not-$(\lambda,\lambda)$-regular ultrafilter. 
\end{corollary}

\begin{proof}
Immediate from Proposition \ref{dsll} and \cite[Theorem 1.7]{C}.
\end{proof} 

In the statements of the next theorems
the word ``productively'', when included within parentheses, can be equivalently inserted or omitted.

\begin{theorem} \label{topprocsing}
For all infinite cardinals $ \lambda  $, $ \mu$, $ \kappa $,
 the following are equivalent:

(i) Every productively $[ \lambda, \mu]$-compact topological space
is (productively) $[ \kappa , \kappa ]$-compact.

(ii) Every productively $[ \lambda, \mu]$-compact family of topological spaces
is productively $[ \kappa , \kappa ]$-compact.

(iii) Every $( \lambda, \mu)$-regular ultrafilter is $( \kappa , \kappa )$-regular.

(iv) Every productively $[ \lambda, \mu]$-compact Hausdorff normal 
topological space with a base of clopen sets
is productively $[ \kappa , \kappa ]$-compact.

(v) Every productively $[ \lambda, \mu]$-compact 
Tychonoff topological group with a base of clopen sets
is (productively) $[ \kappa , \kappa ]$-compact.

If $ \kappa $ is regular, then the preceding conditions are also
equivalent to:

(vi) Every productively $[ \lambda, \mu]$-compact Hausdorff normal 
topological space with a base of clopen sets
is $[ \kappa , \kappa ]$-compact.
\end{theorem}

 \begin{proof} 
Let us denote by (i)$_{\text p}$ Condition (i) when the second occurrence of the word 
``productively'' is included, and simply by (i) when it is omitted. Similarly, for condition (v).

The equivalence of (i)-(iii) has been proved in \cite[Theorem 1]{topproc},
where it has also been proved that,  for $ \kappa $ regular, they are equivalent
to (vi).

Since (ii) $ \Rightarrow $ (i)$_{\mathrm p}$ $ \Rightarrow $ (i)  are trivial,
we get that (i), (ii), (iii), (i)$_{\mathrm p}$ are all equivalent, and equivalent to
(vi) for $\kappa $ regular.

(ii) $\Rightarrow $ (iv) and (ii) $\Rightarrow $ (v)$_{\mathrm p}$ $ \Rightarrow $ (v) are trivial.

If (iii) fails, then there is a
$( \lambda, \mu)$-regular ultrafilter which is not $( \kappa , \kappa )$-regular, thus,
for $\kappa $ singular, 
the space $X$ of Corollary \ref{cor} is 
productively $[ \lambda, \mu]$-compact. 
For $\kappa $ regular, take $X= \kappa $ with the order topology
(see \cite{topproc}).
 $X$ 
is Hausdorff, normal, 
 with a base of clopen sets, but not 
productively $[ \kappa , \kappa ]$-compact, again by Corollary \ref{cor}, 
thus (iv) fails. We have proved (iv) $\Rightarrow $ (iii). 

(v) $\Rightarrow $ (iii) is similar, using 
Corollary \ref{psll}, since  $S_\kappa (\kappa )$ is not
$[\kappa ,\kappa ]$-compact.
\end{proof} 

\begin{theorem} \label{topproc2sing}
For all infinite cardinals $ \lambda  $, $ \mu$, 
and for any family $ (\kappa_i)_{i \in I} $ of infinite cardinals,
 the following are equivalent:

(i) Every productively $[ \lambda, \mu]$-compact topological space
is (productively) 
$[ \kappa_i , \kappa_i ]$-compact for some $i \in I$.

(ii) Every productively $[ \lambda, \mu]$-compact family of topological spaces
is productively 
$[ \kappa_i , \kappa_i ]$-compact for some $i \in I$.

(iii) Every $( \lambda, \mu)$-regular ultrafilter is $( \kappa_i , \kappa_i )$-regular
 for some $i \in I$.

(iv) Every productively $[ \lambda, \mu]$-compact Hausdorff normal 
topological space with a base of clopen sets
is productively 
$[ \kappa_i , \kappa_i ]$-compact for some $i \in I$.

(v) Every productively $[ \lambda, \mu]$-compact 
Tychonoff topological group with a base of clopen sets
is (productively) 
$[ \kappa_i , \kappa_i ]$-compact for some $i \in I$.

If every $ \kappa_i $ is regular, then the preceding conditions are also
equivalent to:

(vi) Every productively $[ \lambda, \mu]$-compact Hausdorff normal 
topological space with a base of clopen sets
is 
$[ \kappa_i , \kappa_i ]$-compact for some $i \in I$.
\end{theorem}

 \begin{proof} 
The equivalence of (i)-(iii) has been proved in 
 \cite[Theorem 3]{topproc}, thus, arguing as in the proof of 
Theorem \ref{topprocsing}, we get that 
(i), (ii), (iii), (i)$_{\mathrm p}$ are all equivalent. 

(ii) $\Rightarrow $ (iv) $\Rightarrow $ (vi) and (ii) $\Rightarrow $ (v)$_{\mathrm p}$ $ \Rightarrow $ (v) are trivial.

If (iii) fails, then there is a
$( \lambda, \mu)$-regular ultrafilter $D$ which for no $i \in I$ is 
$( \kappa_i , \kappa_i )$-regular.
By Proposition \ref{dsll},
for every $i \in I$
the topological space 
$S_{\kappa_i} (\kappa_i )$ is 
 $D$-compact.
Hence $X = \prod_{i \in I} S_{\kappa_i} (\kappa_i )$
is $D$-compact, thus
productively $[ \lambda, \mu]$-compact, by 
\cite[Theorem 1.7]{C}.
However, 
$X$ is a Tychonoff topological group with a base of clopen sets
which for no $i \in I$ is
$[ \kappa_i , \kappa_i ]$-compact,
thus (v) fails.
 We have proved (v) $\Rightarrow $ (iii). 

The proofs of (iv) $\Rightarrow $ (iii) 
and 
(vi) $\Rightarrow $ (iii) 
are similar, using the next proposition.
If (iii) fails, then there is a
$( \lambda, \mu)$-regular ultrafilter $D$ which for no $i \in I$ is 
$( \kappa_i , \kappa_i )$-regular.
By the proof of Theorem \ref{topprocsing},
for every $i \in I$
we have a $D$-compact topological space 
$X_i $ which falsify \ref{topprocsing}(iv),
resp., \ref{topprocsing}(vi). 
Then the space $X = \{x\}  \cupdoty \bigcupdoty _{i \in I} X_i $
we shall construct in the next definition
is $D$-compact, thus
productively $[ \lambda, \mu]$-compact, by 
\cite[Theorem 1.7]{C},
and makes (iv), resp., (vi), fail.
\end{proof} 

\begin{definition}\label{fr}
Given a family $(X_i) _{i \in I} $
of topological spaces, construct their
\emph{Frechet disjoint union}
$X = \{x\}  \cupdoty \bigcupdoty _{i \in I} X_i $
as follows.

Set theoretically, $X$ is the union of
(disjoint copies) of the $X_i$'s, plus a new element
$x$ which belongs to no $X_i$.
The topology on $X$ is the smallest topology which
contains each open set of each  $X_i$,
and which contains 
$\{x\}  \cupdoty \bigcupdoty _{i \in E} X_i $,
for every $E \subseteq I$ such that $I \setminus E$
is finite.
\end{definition}

\begin{proposition}\label{disj}
If $(X_i) _{i \in I} $ is a family of
topological spaces, then
their Frechet disjoint union
$X = \{x\}  \cupdoty \bigcupdoty _{i \in I} X_i $
is $T_0$, $T_1$, Hausdorff, regular, normal, 
$D$-compact (for a given ultrafilter $D$), $[\lambda,\mu]$-compact
(for given infinite cardinals $\lambda $ and $\mu$), 
has a base of clopen sets if and only if  so is (has) each $X_i$.
\end{proposition}

\begin{proof}
Straightforward. We shall comment only on regularity, normality
and $D$-compactness.

For regularity and normality, just observe that 
if $C$ is closed in $X$ and $C$
has nonempty intersection with infinitely many
$X_i$'s, then $x \in C$.

As for $D$-compactness, 
suppose $D$ is over $J$ 
and that each $X_i$ is $D$-compact.
Let $(y_j) _{j \in J} $ be a sequence of elements 
of $X$. If 
$ \{j \in J| y_j \in \bigcup _{i \in F} X_i \} \not\in D $ holds
for every $F \subseteq I$,
 then $(y_j) _{j \in J} $ $D$-converges to $x$.
Otherwise, since $D$ is an ultrafilter, 
hence $\omega$-complete, there exists some
$i \in I$ such that 
$ \{j \in J| y_j \in  X_i \} \in D $.
But then $(y_j) _{j \in J} $ $D$-converges to 
some point of $X_i$, since $X_i$ is supposed to be
$D$-compact. 
\end{proof}

When $ \kappa $ is singular of cofinality $ \omega $,
Condition (vi) in Theorem \ref{topprocsing} is equivalent
to the other conditions.
When each $ \kappa_i $ is either a regular cardinal, or a singular
cardinal of cofinality $ \omega $,
then Condition (vi) in  Theorem \ref{topproc2sing}.
is equivalent
to the other conditions. Proofs shall be given elsewhere.

\end{document}